# ON A PROBLEM WITH CONDITIONS ON ALL BOUNDARY FOR A PSEUDOPARABOLIC EQUATION

## I.G.Mamedov

*A.I.Huseynov Institute of Cybernetics of NAS of Azerbaijan. Az 1141, Azerbaijan, Baku st. B. Vahabzade, 9*
*E-mail: ilgar-mammadov@rambler.ru*

### Abstract

*A problem with non-classical conditions on all the boundary not requiring agreement conditions is considered for a sixth order pseudoparapolic equation. The equivalence of these conditions with the classic boundary condition is substantiated in the case if the solution of the stated problem is sought in S.L.Sobolev isotropic space $W_p^{(3,3)}(G)$.*

## Problem statement

Consider the equation

$$(V_{3,3}u)(x,y) \equiv \sum_{i,j=0}^{3} a_{i,j}(x,y) D_x^i D_y^i u(x,y) = Z_{3,3}(x,y) \in L_p(G), \quad (1)$$

where $a_{3,3}(x,y) \equiv 1$.

Here $u(x,y)$ is a desired function defined on $G$; $a_{i,j}(x,y)$ are the given measurable function on $G = G_1 \times G_2$, where $G_k = (0, h_k)$, $k = 1,2$; $Z_{3,3}(x,y)$ is a given measurable function on $G$; $D_t^k = \partial^k / \partial t^k$ is a generalized differentiation operator in S.L.Sobolev sense $D_t^0$ is an identity transformation operator.

Equation (1) is a hyperbolic equation possessing two real characteristics $x = const$, $y = const$, the first and second of which is three-fold. Therefore, in some sense we can consider equation (1) as a pseudoparabolic equation [1]. This equation is a generalization of the equation of thin spherical shell bending [2. p. 258].

In this paper we consider equation (1) in the general case when the coefficients $a_{i,j}(x,y)$ are non-smooth functions satisfying only the conditions:

$$a_{i,j}(x,y) \in L_p(G), \ i=\overline{0,2}, \ j=\overline{0,2};$$

$$a_{3,j}(x,y) \in L_{\infty,p}^{x,y}(G), \ j=\overline{0,2}.$$

$$a_{i,3}(x,y) \in L_{p,\infty}^{x,y}(G), \ i=\overline{0,2}.$$

Therewith the important principal moment is that the equation under consideretion possesses discontinuous coefficients that satisfy only some $p$-integrability and boundedness conditions, i.e. the considered differential operator $V_{3,3}$ has no traditional conjugated operator.

Under this conditions, the solution $u(x,y)$ of equation (1) will be sought in S.L.Sobolev space

$$W_p^{(3,3)}(G) = \{u(x,y) : D_x^i D_y^j u(x,y) \in L_p(G), i=\overline{0,3}, j=\overline{0,3}\},$$

where $1 \le p \le \infty$. We'll define the norm in the isotropic space $W_p^{(3,3)}(G)$ by the equality

$$\|u\|_{W_p^{(3,3)}(G)} = \sum_{i,j=0}^{3} \|D_x^i D_y^j u\|_{L_p(G)}.$$

We can give the classic form conditions on all the boundary for equation (1) as follows [3]:

$$\begin{cases} u(0,y) = \varphi_1(y); \ u(x,0) = \psi_1(x); \\ u(h_1,y) = \varphi_2(y); u(x,h_2) = \psi_2(x); \\ \left.\dfrac{\partial u(x,y)}{\partial x}\right|_{x=0} = \varphi_3(y); \left.\dfrac{\partial u(x,y)}{\partial y}\right|_{y=0} = \psi_3(x); \end{cases} \quad (2)$$

where $\varphi_k(y), \psi_k(x), k=\overline{1,3}$ are the given measurable functions on $G$. Obviously, in the case of conditions (2) in addition to conditions

$$\varphi_k(y) \in W_p^{(3)}(G_2) \equiv \{\widetilde{\varphi}(y) : D_y^j \widetilde{\varphi}(y) \in L_p(G_2), j=\overline{0,3}\}, \quad 1 \le p \le \infty;$$

$$\psi_k(x) \in W_p^{(3)}(G_1) \equiv \{\widetilde{\psi}(x) : D_x^j \widetilde{\psi}(x) \in L_p(G_1), j=\overline{0,3}\}, \quad 1 \le p \le \infty,$$

the given functions satisfy also the following agreement conditions (see. Fig. 1):

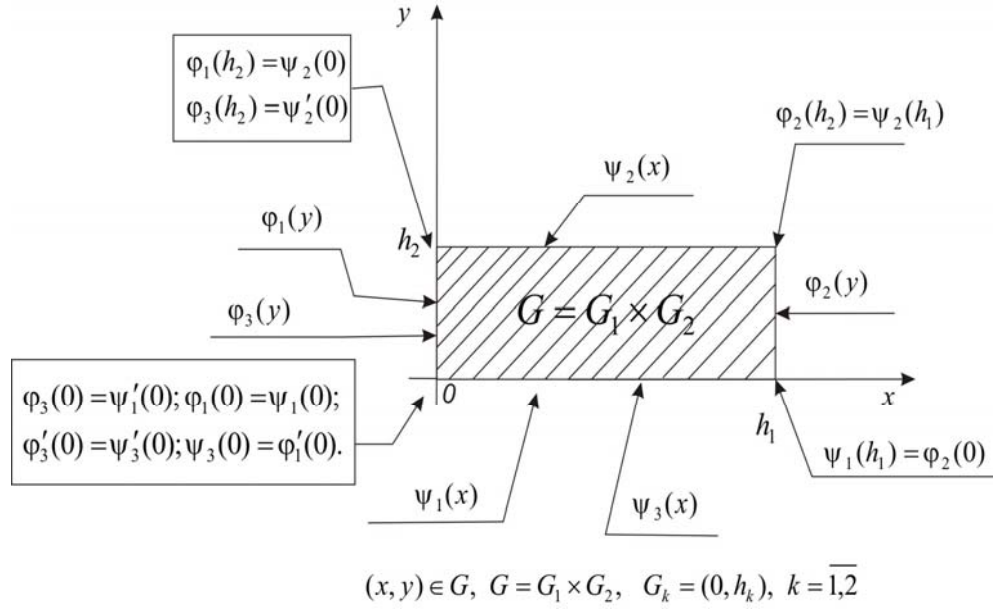

$(x, y) \in G, \ G = G_1 \times G_2, \ G_k = (0, h_k), \ k = \overline{1,2}$

**Fig.1. Geometric interpretation of classical conditions on all the boundary.**

$$\begin{cases} \varphi_1(h_2) = \psi_2(0); \ \varphi_3(0) = \psi_1'(0); \\ \varphi_1(0) = \psi_1(0); \ \varphi_3(h_2) = \psi_2'(0); \\ \psi_1(h_1) = \varphi_2(0); \ \varphi_3'(0) = \psi_3'(0); \\ \varphi_2(h_2) = \psi_2(h_1); \ \psi_3(0) = \varphi_1'(0); \ \psi_3(h_1) = \varphi_2'(0) \end{cases} \quad (3)$$

Obviously, conditions (2) are close to boundary conditions of the Dirichlet problem from [4-7].

Consider the following non-classical boundary conditions:

$$\begin{cases} V_{0,j}u \equiv D_y^j u(0,0) = Z_{0,j} \in R, \ j = \overline{0,2}; \\ (V_{0,3}u)(y) \equiv D_y^3 u(0, y) = Z_{0,3}(y) \in L_p(G_2); \\ V_{1,j}u \equiv D_x D_y^j u(0,0) = Z_{1,j} \in R, \ j = \overline{0,2}; \\ (V_{1,3}u)(y) \equiv D_x D_y^3 u(0, y) = Z_{1,3}(y) \in L_p(G_2); \\ V_{2,j}u \equiv D_x^2 D_y^j u(0,0) = Z_{2,j} \in R, \ j = \overline{0,1}; \\ (V_{3,j}u)(x) \equiv D_x^3 D_y^j u(x,0) = Z_{3,j}(x) \in L_p(G_1), \ j = \overline{0,1}; \\ V_{0,j}^{(h_1)}u \equiv D_y^j u(h_1,0) = Z_{0,j}^{(h_1)} \in R, \ j = \overline{0,2}; \\ (V_{0,3}^{(h_1)}u)(y) \equiv D_y^3 u(h_1, y) = Z_{0,3}^{(h_1)}(y) \in L_p(G_2); \\ V_{i,0}^{(h_2)}u \equiv D_x^i u(0,h_2) = Z_{i,0}^{(h_2)} \in R, i = \overline{0,2}; \\ (V_{3,0}^{(h_2)}u)(x) \equiv D_x^3 u(x,h_2) = Z_{3,0}^{(h_2)}(x) \in L_p(G_1). \end{cases} \quad (4)$$

If the function $u(x,y) \in W_p^{(3,3)}(G)$ is a solution of the classical form problem (1), (2), with conditions on all the boundary then it is also a solution of problem (1),(4) for $Z_{i,j}$, $Z_{i,j}^{(h_1)}$, $Z_{i,j}^{(h_2)}$, determined by the following equalities:

$$Z_{0,0} = \varphi_1(0) = \psi_1(0); \quad Z_{0,1} = \varphi_1'(0) = \psi_3(0); \quad Z_{0,2} = \varphi_1''(0);$$

$$Z_{0,3}(y) = \varphi_1'''(y); Z_{1,0} = \varphi_3(0) = \psi_1'(0); \quad Z_{1,1} = \varphi_3'(0) = \psi_3'(0); \quad Z_{1,2} = \varphi_3''(0);$$

$$Z_{1,3}(y) = \varphi_3'''(y); \quad Z_{2,0} = \psi_1''(0); \quad Z_{3,0}(x) = \psi_1'''(x); \quad Z_{2,1} = \psi_3''(0); \quad Z_{3,1}(x) = \psi_3'''(x);$$

$$Z_{0,0}^{(h_1)} = \varphi_2(0) = \psi_1(h_1); \quad Z_{0,1}^{(h_1)} = \varphi_2'(0) = \psi_3(h_1); \quad Z_{0,2}^{(h_1)} = \varphi_2''(0); \quad Z_{0,3}^{(h_1)}(y) = \varphi_2'''(y);$$

$$Z_{0,0}^{(h_2)} = \psi_2(0) = \varphi_1(h_2); \quad Z_{1,0}^{(h_2)} = \psi_2'(0) = \varphi_3(h_2); \quad Z_{2,0}^{(h_2)} = \psi_2''(0); \quad Z_{3,0}^{(h_2)}(x) = \psi_2'''(x).$$

It is easy to prove that the inverse one is also true. In other words, if the function $u \in W_p^{(3,3)}(G)$ is a solution of problem (1), (4) (see. Fig. 2), then it is also a solution of problem (1), (2) for the following functions:

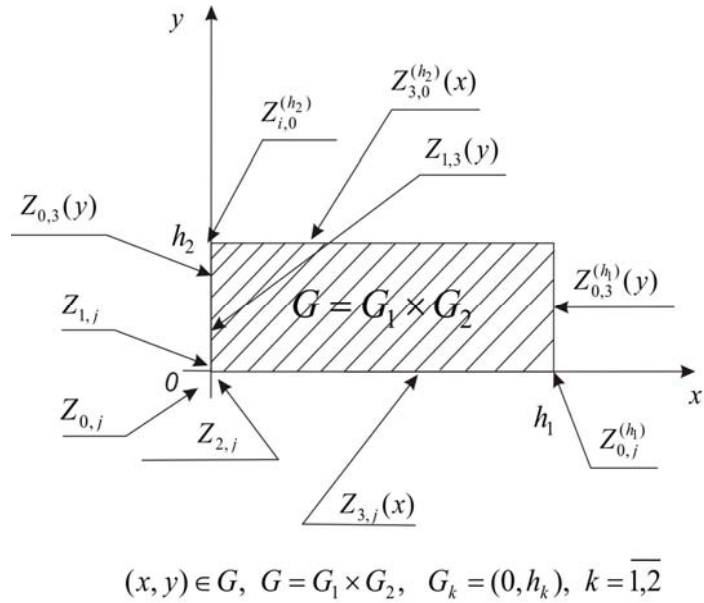

$(x,y) \in G$, $G = G_1 \times G_2$, $G_k = (0, h_k)$, $k = \overline{1,2}$

**Fig.2. Geometrical interpretation of conditions on all the boundary in non classic treatment.**

$$\varphi_1(y) = Z_{0,0} + yZ_{0,1} + \frac{y^2}{2}Z_{0,2} + \frac{1}{2}\int_0^y (y-\tau)^2 Z_{0,3}(\tau)d\tau; \qquad (5)$$

$$\varphi_2(y) = Z_{0,0}^{(h_1)} + yZ_{0,1}^{(h_1)} + \frac{y^2}{2}Z_{0,2}^{(h_1)} + \frac{1}{2}\int_0^y (y-\xi)^2 Z_{0,3}^{(h_1)}(\xi)d\xi; \qquad (6)$$

$$\varphi_3(y) = Z_{1,0} + yZ_{1,1} + \frac{y^2}{2}Z_{1,2} + \frac{1}{2}\int_0^y (y-\eta)^2 Z_{1,3}(\eta)d\eta; \quad (7)$$

$$\psi_1(x) = Z_{0,0} + xZ_{1,0} + \frac{x^2}{2}Z_{2,0} + \frac{1}{2}\int_0^x (x-\tau)^2 Z_{3,0}(\tau)d\tau; \quad (8)$$

$$\psi_2(x) = Z_{0,0}^{(h_2)} + xZ_{1,0}^{(h_2)} + \frac{x^2}{2}Z_{2,0}^{(h_2)} + \frac{1}{2}\int_0^x (x-\nu)^2 Z_{3,0}^{(h_2)}(\nu)d\nu; \quad (9)$$

$$\psi_3(x) = Z_{0,1} + xZ_{1,1} + \frac{x^2}{2}Z_{2,1} + \frac{1}{2}\int_0^x (x-\mu)^2 Z_{3,1}(\mu)d\mu; \quad (10)$$

Note that the functions (5)-(10) possess one important property, more exactly, agreement condition (3) for all $Z_{i,j}$, $Z_{i,j}^{(h_1)}$, $Z_{i,j}^{(h_2)}$ having the above-stated properties are fulfilled for them automatically. Therefore equalities (5)-(10) may be considered as a general form of all the functions $\varphi_k(y)$, $\psi_k(x), k = \overline{1,3}$, satisfying agreement conditions (3).

So, classic type problems (1), (2) and of the form (1), (4) with conditions on all the boundary are equivalent in general case. However, the non - classic problem (1), (4) with conditions on the boundary is more natural by statement than problem (1), (2).

This is connected with the fact that in the statement of problem (1), (4) with conditions on all the boundary, the right parts of boundary conditions don't require additional conditions of agreement type. Note that various boundary value problems with non classic conditions not requiring agreement conditions were substantiated in the author's papers [8-9].